\documentclass[journal]{IEEEtran}

\usepackage[english]{babel}
\usepackage{ifpdf}

\ifCLASSINFOpdf
 \usepackage[pdftex]{graphicx}
\else
\usepackage[dvips]{graphicx}
\fi

\usepackage{amsmath}
\usepackage{amssymb}
\addto\captionsenglish{}
\usepackage{url}
\usepackage[utf8]{inputenc}
\usepackage{algpseudocode}
\usepackage{algorithm}
\usepackage{cite}
\usepackage{enumitem}

\usepackage{eurosym}
\usepackage{amstext} 
\DeclareRobustCommand{\officialeuro}{%
  \ifmmode\expandafter\text\fi
  {\fontencoding{U}\fontfamily{eurosym}\selectfont e}}

\begin{document}

\title{Coalitional game based cost optimization of energy portfolio in smart grid communities}
\author{Adriana~Chi\c{s},~\IEEEmembership{Student Member,~IEEE,}
		and Visa~Koivunen,~\IEEEmembership{Fellow,~IEEE}
\thanks{The authors are with the Department of Signal Processing and Acoustics, Aalto University, FI-02150 Espoo, Finland (e-mail: adriana.chis@aalto.fi, visa.koivunen@aalto.fi)}}

\date{}

\maketitle

\begin{abstract}
In this paper we propose two novel coalitional game theory based optimization methods for minimizing the cost of electricity consumed by households from a smart community. Some households in the community may own renewable energy systems (RESs) conjoined with energy storing systems (ESSs). Some other residences own ESSs only, while the remaining households are simple energy consumers. We first propose a coalitional cost optimization method in which RESs and ESSs owners exchange energy and share their renewable energy and storage spaces. We show that by participating in the proposed game these households may considerably reduce their costs in comparison to performing individual cost optimization. We further propose another coalitional optimization model in which RESs and ESSs owning households not only share their resources, but also sell energy to simple energy consuming households. We show that through this energy trade the RESs and ESSs owners can further reduce their costs, while the simple energy consumers also gain cost savings. The monetary revenues gained by the coalition are distributed among its members according to the Shapley value. Simulation examples show that the proposed coalitional optimization methods may reduce the electricity costs for the RESs and ESSs owning households by 20\%, while the sole energy consumers may reduce their costs by 5\%. 
\end{abstract}


\IEEEpeerreviewmaketitle

\section{Introduction}
\label{section:1}

%
%
%

\IEEEPARstart{O}{ptimization} of energy consumption imposes changes both at the energy supply side and at the energy demand side. One key global challenge existing today is the achievement of sustainable energy production. More and more energy users choose to install renewable energy systems (RESs) which seem to be a sustainable solution for the current environmental problems, but also for other troubling issues existing in the power network today. Installation of RESs at distribution level can result in significant reduction of power transmission losses, lower operational costs and an overall cut of electricity costs. 
 	
Energy consumption optimization at the distribution level can also be obtained through demand side management (DSM) methods. DSM refers to the modification of the energy consumption patterns of the end-users with the purpose of lowering costs, reducing load peaks on the grid and increasing grid reliability. Energy storing systems (ESSs) together with a smart metering infrastructure can be an easy and efficient way of implementing such methods without the need of highly modifying the energy consumption patterns of energy consumers such as residential households and business centers.   

One major disadvantage of the RESs is the intermittent nature of their energy production, i.e. not providing  reliable energy production. ESSs can be used in conjunction with RESs in order to store the surplus of the produced energy, but also for DSM operations. However, installation of ESSs with high capacities can be very costly. One way of addressing these drawbacks of RESs and ESSs is by employing cooperation among RESs and ESSs owners. Cooperation can improve the integration of renewable energy, reduce the need of large ESSs and provide cost savings for the participants in the cooperation.

Cooperative methods for energy trading have been studied before. The studies in \cite{A1_2,B2_2,C3_2,D4_2,E5_2,F6_2,G7_2} approach the problem of energy trading and cooperation between microgrids. Methods based on Nash bargaining theory to incentivize the collaboration between microgrids are proposed in \cite{A1_2} and \cite{B2_2}. An approach for cooperative power control in a network of microgrids is proposed in \cite{D4_2} with the purpose of optimizing the use of renewable energy for meeting the demand and to enhance the reliability of the network. In \cite{F6_2} prospect theory has been proposed to formulate a static game between autonomous microgrids that are trading energy.  

In this paper we propose two novel coalitional game theoretic based models for optimizing energy portfolios within a smart grid community by sharing and trading energy among households. We consider a community of households some of which own RESs together with ESSs. Some other households in the community own ESSs only, while the remaining households are simple energy consumers. We assume that the households are equipped with smart energy management meters that can predict with sufficient accuracy their energy demand profiles and the profiles of renewable energy production during a finite time period ahead. The households are connected to a centralized energy management unit that finds the solution of the proposed energy cost optimization problems and controls the flow of energy within the community through a two-way communication system. In case of insufficient renewable energy resources, the electricity demands of the households are fulfilled with electricity bought from the utility company. Our main contributions to this paper are described below: 
\begin{itemize} [leftmargin=*] 
\item We formulate three cost optimization problems: 1) a method through which households that own RESs and ESSs can individually optimize their costs; 2) a coalitional optimization method in which the RESs and ESSs owners may exchange energy and share the produced renewable energy and storage spaces; 3) a coalitional optimization method in which the RESs and ESSs may exchange energy, share their resources and storage spaces and also sell energy to sole energy consuming households in the community. The sole energy consumers are not participants in the coalitional game. They only act as agents that buy energy from the coalitional group of RESs and ESSs owners at a cheaper price than the one offered by the utility company. The proposed optimization problems are convex and linear and can be easily solved through standard linear programming methods. 
\item We model a coalitional game among the RESs and ESSs owners. The revenue gained by the coalitional optimization is divided among the coalition members in a fair manner, according to the amount of contribution that each member brings to the coalition. A method based on Shapley value is used to calculate the payoff earned by each member. 
\item We perform extended simulations and show that through the collaborative method the smart households owning RESs and ESSs can obtain significant financial benefits in comparison with the individual cost optimization. We also show that by selling energy to the sole energy consumers at a smaller price than the one offered by the utility company, the RESs and ESSs owners can further reduce their electricity costs, while the sole energy consumers also get financial benefits.
\end{itemize}

Game theoretic methods for energy trading in smart grids have also been proposed in \cite{H8_2,L12_2,I9_2,K11_2,J10_2}. A noncooperative game among storage units has been proposed in \cite{H8_2}. A cooperative game between households has been proposed in \cite{L12_2} with the scope of flattening the community's load.\cite {I9_2} proposes a noncooperative and a cooperative game among energy users that own distributed energy generation and storages. A coalition-formation game is proposed in \cite{K11_2} for microgrids to trade power in order to reduce the load and the power losses over the power grid. \cite{J10_2} proposes a coalitional game in which microgrids supply a energy to a shared facility controller with the purpose of gaining some revenues.

The coaliotional method proposed in this paper has the following novel aspects: in the proposed game the energy exchange among the coalition members is completely free, the members share their renewable resources and storage spaces and perform DSM services for each other; for further reducing the costs, sole energy consuming households are introduced in the optimization; the cost savings obtained by the sole energy consumers result from the cost minimization achieved by the RESs and ESSs owners through selling them energy at a price lower than the one offered by the utility company. 

Preliminary results related to this work were presented in \cite{M13_2}. This paper extends these results by introducing the sole energy consuming households to the coalitional optimization. The optimization problems are improved by adding operational costs and penalty terms. We perform extended simulations to demonstrate the performance of the proposed methods.

The simulations show that over a long term period the coalitional game among the RESs and ESSs owning households can result in a cost reduction of approximately 13.5\% more in comparison with individual cost optimization. The introduction of the sole energy consuming households can achieve a further cost reduction of 6.5\% for the RESs and ESSs owners, resulting in a total 20\% reduction of the cost for these residences. The sole energy consumers can also achieve a cost reduction of about 5\% in comparison to buying all needed electricity from the utility company. 

The rest of this paper is organized as follows: Section II presents the system model, Section III presents the three cost optimization methods. The coalitional game model for the RESs and ESSs owning households and the Shapley value method are presented in Section IV, Section V shows the simulation results and Section VI presents the conclusions.

\section{System model}
\label{section:2}

We consider a community of $N$ households that can trade and exchange energy with each other through a central energy control unit. We denote by $\mathcal{N}$ the set of households in the community and by $n$ we denote the index of a household from this set.  A subset of smart residences from the community, denoted by $\mathcal{M}$, $\mathcal{M}\subseteq \mathcal{N}$, with cardinality $|\mathcal{M}|=M$, either produce renewable energy,  which implies owning an ESS as well, or own ESSs only. We denote the index of a household from this set by $m$. We denote by $\mathcal{P}$, with cardinality $|\mathcal{P}|=P$, $\mathcal{P}=\mathcal{N} \setminus \mathcal{M}$, the remaining subset of households from the community which are sole energy consumers. We denote the index of a household from this set by $p$. The energy optimization is performed over a finite time horizon $\mathcal{T}$ of length $ T$, which is divided into equally long time-slots denoted by $t$, $\mathcal{T}=[t,t=1,\ldots,T]$. 

\subsection{The Energy Model}    	

The electricity demand of each household in the community is considered to be known ahead for the entire period $\mathcal{T}$ and it is not flexible. The set of electricity demands of one residence in the community is denoted by \emph{$\mathbf{u}_n=[u_n(t),t=1,\ldots,T],n\in\mathcal{N}$}. The set of per-time-slot amounts of renewable energy, \emph{$\mathbf{w}_m=[w_m(t),t=1,\ldots,T]$}, produced by those households $m \in \mathcal{M}$ owning RESs is considered predicted with sufficient accuracy over the whole period $\mathcal{T}$. Note that for those households from the set $\mathcal{M}$ that do not own RESs, but only ESSs, the renewable energy vector is zero, i.e. $\mathbf{w}_m=0$. In this paper, the problem of exchanging and trading the energy among the households is formulated as a coalitional game in which the participants may produce, store and exchange energy. Hence, we denote by \emph{$\mathbf{a}_n=[a_n(t),t=1,\ldots,T]$} the set of energy amounts that a household $n$ from the community may exchange with the rest of households from the community within period $\mathcal{T}$. This variable may have either positive or negative values. If  $a_n(t)>0$ in a certain time slot \emph{t}, it means that household $n$ provides this amount of energy to the rest of the households in the community, while if $a_n(t)<0$ then household $n$ receives this amount of energy from the other members of the community. The set of energy amounts that a household $n$ may have to purchase from the utility company during period $\mathcal{T}$ is denoted by \emph{$\mathbf{b}_n=[b_n(t),t=1,\ldots,T]$}. The energy purchased by a household from the utility company in a time-slot $t$ must obey the energy consumption constraint. This constraint is specific for each type of optimization and for each set of households in the community.

The households belonging to set $\mathcal{M}$, owning RESs and/or ESSs, can optimize their cost individually, without collaborating with other members of the community. In this case the following energy consumption constraint may be imposed:
\begin{equation}\label{eq:eq1}
u_m(t)\hspace{-0.5mm}-\hspace{-0.5mm}w_m(t)\hspace{-0.5mm}-\hspace{-0.5mm}b_m(t)\hspace{-0.5mm}+\hspace{-0.5mm}r_m(t)\leq 0,\hspace{1mm}\forall t\in\mathcal{T},\forall\hspace{1mm} m\in \mathcal{M}.
\end{equation}
This inequation states that in a time-slot $t$ the electricity demand of a household, $u_m(t$), must be fulfilled by the available renewable energy, $w_m(t)$, the energy purchased from the grid, $b_m(t)$, and by energy from in the storage, $r_m(t)$.

In the collaborative case, the households belonging to set $\mathcal{M}$ may obey the following energy consumption constraint:
\begin{equation}\label{eq:eq2}
u_m(t)\hspace{-0.5mm}-\hspace{-0.5mm}w_m(t)\hspace{-0.5mm}-\hspace{-0.5mm}b_m(t)\hspace{-0.5mm}+\hspace{-0.5mm}r_m(t)+\hspace{-0.5mm}a_m(t)\hspace{-0.5mm}\leq\hspace{-0.5mm}0,\forall t\hspace{-0.5mm}\in\hspace{-0.5mm}\mathcal{T},\forall m\hspace{-0.5mm}\in\mathcal{M}.
\end{equation}
The interpretation of this constraint is similar to that of constraint (\ref{eq:eq1}). The difference is that in the collaborative scenario one participant to the optimization may also receive or transfer an amount of electricity $a_m(t)$ from or to other households in the community.

The energy consumption constraint for the pure consuming households $\mathcal{P}$ is defined as follows:
\begin{equation}\label{eq:eq3}
u_p(t)\hspace{-0.5mm}-\hspace{-0.5mm}b_p(t)\hspace{-0.5mm}+\hspace{-0.5mm}a_p(t)= 0,\hspace{1mm}\forall t\in\mathcal{T}, \forall \hspace{1mm} p\in \mathcal{P}.
\end{equation}
This category of households can only receive energy from the other members of the community. This imposes the following constraint:
\vspace{-1.5mm}
\begin{equation}\label{eq:eq4}
a_p(t)\leq 0,\hspace{1mm}\forall t\in\mathcal{T}, \forall \hspace{1mm} p\in \mathcal{P}.
\end{equation}
Alternatively, they can buy it from the utility company. For all the households in the community the amount of energy purchased from the power grid in time-slot $t$, $b_n(t)$, may have only positive or zero values:
\begin{equation}\label{eq:eq5}
b_n(t)\geq 0,\hspace{1mm}\forall t\in\mathcal{T}, \forall \hspace{1mm} n\in \mathcal{N}.
\end{equation}
In the proposed model the households cannot sell back electricity to the power grid.

\subsection{The Energy Storage Model} 

All households from subset $\mathcal{M}$ own ESSs. We denote by $C_m$ be the maximum storing capacity associated with each $ESS_m$. The maximum amount of energy that can be charged or discharged from a storage is limited by the charging/discharging rate. In this work the charging rate is equal to the discharging rate and we denote this parameter by $\rho_m$. Another parameter that characterizes a storage system is the storage leakage factor, $\eta_m$, which shows the proportion of energy that a storage loses during a time unit. This parameter has values between 0 and 1, but typically $\eta_m\ll 1$. 

We further define by \emph{$\mathbf{r}_m=[r_m(t),t=1,\ldots,T]$} the set of amounts of energy charged or discharged from a storage unit $ESS_m$ during period $\mathcal{T}$. This value can have either positive or negative values. If $r_m(t)>0$ then energy is charged into the storage, while if $r_m(t)<0$ then energy is being discharged from the storage in time-slot $t$. The amount of energy charged or discharged from storage is bounded by the charging/discharging rate $\rho_m$:
\begin{equation}\label{eq:eq6}
-\rho_m \leq r_m(t)\leq \rho_m,\hspace{1mm}\forall t\in\mathcal{T}, \forall \hspace{1mm} m\in \mathcal{M}.
\end{equation}
Let \emph{$\mathbf{s}_m=[s_m(t),t=1,\ldots,T]$} be the energy storage vector containing the total amounts of energy stored in an $ESS_m$ at the end of each time slot. The dynamics of the $ESS_m$ is defined by the following equation:
\begin{equation}\label{eq:eq7}
s_m(t)=(1-\eta_m)s_m(t-1)+r_m(t),\hspace{1mm}\forall t\in\mathcal{T}, \forall \hspace{1mm} m\in \mathcal{M},
\end{equation}
where $s_m(0)$ would be the initial storage value, i.e the amount of energy remained in storage at the end of the previous optimization period. The amount of energy existing in storage at any time-slot must obey the storage capacity constraint:
\begin{equation}\label{eq:eq8}
0\leq s_m(t)\leq C_m,\hspace{1mm}\forall t\in\mathcal{T}, \forall \hspace{1mm} m\in \mathcal{M}.
\end{equation}

\section{The cost minimization problem}
\label{section:3}

We propose two coalitional game based approaches for minimizing energy costs among cooperating households. In the first proposed coalitional optimization problem the households that own RESs and ESSs share their energy storage units and their renewable energy resources and exchange energy among themselves. We show that through participating to the proposed coalitional optimization the households that own RESs and/or ESSs may significantly reduce their electricity costs in comparison to the case in which they would individually optimize their energy costs, using only their own renewable resources and energy storages. Moreover, we propose another coalitional game model in which we show that by including the households that only consume energy in the optimization model, the households owning RESs and ESSs may reduce their costs even further by selling energy to these households at a price lower than the one offered by the utility company. By participating in this optimization, the sole energy consuming households can also reduce their energy costs.

The market electricity prices per unit of energy are given by the utility company and they are known ahead for each time-slot $t$ in the period $\mathcal{T}$. We denote this set of electricity prices by \emph{$\boldsymbol{\xi}=[\xi(t),t=1,\ldots,T]$}.

In the following subsections we describe the proposed cost minimization problems. First we formulate a cost minimization problem which can be performed by each household $m\in\mathcal{M}$ individually, without collaborating. Then we formulate  the coalitional cost minimization problem that includes only the residences that own RESs and ESSs. Finally, we formulate the coalitional cost minimization problem which includes the households that are sole electricity consumers as well.

\subsection{Individual cost minimization problem}
 We define the cost of electricity bought from the utility company by any household $n\in \mathcal{N}$ as: 
\vspace{-2mm}
\begin{equation}\label{eq:eq9}
C_n^{\text{grid}}=\sum_{t=1}^T \xi(t)b_n(t),\hspace{1mm} \forall \hspace{1mm} n\in \mathcal{N}.
\end{equation}

In addition to this cost, the households $m\in\mathcal{M}$ are also incurred with an energy storage cost. It is well known that while being charged and discharged multiple times, storage units suffer a certain degradation and the lifetime of a storage unit is limited to a number of charging/discharging cycles. We define by $\pi$ the storage degradation price corresponding to charging/ischarging one unit of energy. The resulting $ESS_m$ storage degradation cost within period $\emph{T}$ can be defined as:
\vspace{-2mm}
\begin{equation}\label{eq:eq10}
C_m^{\text{storage}}=\pi\sum_{t=1}^T |r_m(t)|,\hspace{1mm} \forall \hspace{1mm} m\in \mathcal{M}.
\end{equation}

Households $m\in \mathcal{M}$, owning RESs and/or ESSs can individually optimize their costs, without collaborating with the other members of the community. We denote by $c_m^{\text{indiv}}$ the individual cost of a household from set $\mathcal{M}$:
\begin{equation}\label{eq:eq21}
 c_m^{\text{indiv}}=C_m^{\text{grid}}+C_m^{\text{storage}}.
\end{equation}

The individual cost minimization problem for one single household $m\in\mathcal{M}$ can be stated as follows:

\vspace{2mm}
\hspace{-3mm}$\min\limits_{\mathbf{b}_m, \mathbf{r}_m, \mathbf{s}_m} c_m^{\text{indiv}} +\sigma \sum\limits_{t=1}^T\left[b_m(t)+w_m(t)\hspace{-1mm}-\hspace{-1mm}u_m(t)\hspace{-1mm}-\hspace{-1mm}r_m(t)\right]$,\\

\hspace{-3.5mm}such that the constraints (\ref{eq:eq1}), (\ref{eq:eq5})-(\ref{eq:eq8}) are satisfied. The last term of the objective function: $\sigma \sum\limits_{t=1}^T\left[b_m(t)\hspace{-1mm}+\hspace{-1mm}w_m(t)\hspace{-1mm}-\hspace{-1mm}u_m(t)\hspace{-1mm}-\hspace{-1mm}r_m(t)\right]$, is a penalty term which is closely related to the energy consumption constraint (\ref{eq:eq1}). The purpose of this penalty is to make sure that the renewable energy produced in a time-slot is either consumed or stored in the $ESS_m$. By $\sigma$ we define a penalty price for a unit of renewable energy that may be wasted. As mentioned in Section II, the set $\mathcal{M}$ of households contains households that own both RESs and ESSs and households that own ESSs only. For the second category $\mathbf{w}_m=0$ and the optimization is performed by using their energy storages alone. 

The set of solution variables of the individual cost minimization problem is $\{\mathbf{b}_m, \mathbf{r}_m, \mathbf{s}_m\}$.

\subsection{Coalitional cost minimization problem for households owning RESs and/or ESSs}

In this scenario the households belonging to the set $\mathcal{M}$ collaborate and share their renewable resources and storage units in order to jointly optimize their costs. Besides the cost of electricity purchased from the power grid and the cost incurred for storage degradation, a third type of costs shall be added to the optimization problem in the collaborative scenario. This cost represents the cost for operating the central control unit. The optimization is performed by a central energy management unit that performs the optimization and controls the energy flow and transfer within the community.  For this, we define by $\tau$ the price charged for transferring one unit of energy from one household to another. The overall cost incurred by a household $m$ during the period $\mathcal{T}$ for the energy transfer operations is equal to:
\vspace{-2mm}
\begin{equation}\label{eq:eq11}
C_m^{\text{operation}}=\tau\sum_{t=1}^T |a_m(t)|,\hspace{1mm} \forall \hspace{1mm} m\in \mathcal{M}.
\end{equation}

We denote by $c_{\mathcal{M}}^{\text{ResEss}}$ the overall cost jointly incurred by all households owning RESs and ESSs: 

\vspace{-5mm}
\begin{equation}\label{eq:eq22}
c_{\mathcal{M}}^{\text{ResEss}}=\sum\limits_{m=1}^M C_m^{\text{grid}}\hspace{-1mm}+\hspace{-1mm}\sum\limits_{m=1}^M C_m^{\text{storage}}\hspace{-1mm}+\hspace{-1mm}\sum\limits_{m=1}^M C_m^{\text{operation}}. 
\end{equation}

The coalitional cost minimization problem for the households $m\in\mathcal{M}$ that own RESs and/or ESSs is formulated as:

\vspace{2mm}
\hspace{-3.5mm}$\min\limits_{\{\mathbf{b}_m, \mathbf{r}_m, \mathbf{s}_m, \mathbf{a}_m \}_{m=1}^M} \hspace{0.5mm} c_{\mathcal{M}}^{\text{ResEss}}+$\\

\hspace{19mm}$\sigma\hspace{-1.5mm}\sum\limits_{m=1}^M\hspace{-0.5mm}\sum\limits_{t=1}^T[b_m(t)\hspace{-0.5mm}+\hspace{-0.5mm}w_m(t)-\hspace{-0.5mm}u_m(t)\hspace{-0.5mm}-\hspace{-0.5mm}a_m(t)\hspace{-0.5mm}-\hspace{-0.5mm}r_m(t)]$, \\

\hspace{-3.5mm}such that the constraints (\ref{eq:eq2}), (\ref{eq:eq5})-(\ref{eq:eq8}) are satisfied. A penalty term is added to the objective also in case of the coalitional optimization. The penalty ensures that the produced renewable energy is either consumed, either transferred to the storage. In this case penalty term is expressed by: $\sigma \sum\limits_{m=1}^M  \sum\limits_{t=1}^T\left[b_m(t)\hspace{-1mm}+\hspace{-1mm}w_m(t)\hspace{-1mm}-\hspace{-1mm}u_m(t)\hspace{-1mm}-\hspace{-1mm}a_m(t)-\hspace{-1mm}r_m(t)\right]$ and it works in corroboration with energy consumption constraint (\ref{eq:eq2}).

To ensure a correct functionality of the optimization problem we add another balance constraint which makes sure that the total amount of energy that is given away by some households is equal to the total amount of energy received by the rest of households from the set $\mathcal{M}$:
\vspace{-2mm}
\begin{equation}\label{eq:eq14} 
\sum_{m=1}^M a_m(t)=0. 
\end{equation}

The final set of constraints for this problem is (\ref{eq:eq2}), (\ref{eq:eq5})-(\ref{eq:eq8}), (\ref{eq:eq14}). The set of solution variables is $\{\mathbf{b}_m, \mathbf{r}_m, \mathbf{s}_m, \mathbf{a}_m \}_{m=1}^M$.

\subsection{Coalitional cost minimization problem for the whole community of households}

In this scenario we include as well the households $p\in\mathcal{P}$ which are sole energy consumers. This set of households participates to the community coalitional energy optimization by just buying energy from those owning RESs and ESSs. 

We define by \emph{$\boldsymbol{\lambda}=[\lambda(t),t=1,\ldots,T]$} the set of prices corresponding to the energy sold by the households $m\in \mathcal{M}$, owning RESs and ESSs, to the households $p\in \mathcal{P}$. These prices are lower than the per-time-slot prices offered by the utility prices: $\lambda(t)=\alpha\xi(t)$, where $0\leq\alpha\leq 1$. For a household $p$ the cost of purchasing energy from the community can be stated as:
\vspace{-3mm}
\begin{equation}\label{eq:eq16} 
C_p^{\text{purchase}}=\sum_{t=1}^T\lambda(t)[-a_p(t)],\hspace{1mm} \forall \hspace{1mm} p\in \mathcal{P}.
\end{equation}

We denote by $c_{\mathcal{M}}^{\text{community}}$ the overall cost incurred by those households from the community owning RESs and ESSs after selling energy to those that are sole energy consumers. This cost is expressed as: 
\vspace{-2mm}
\begin{equation}\label{eq:eq23} 
c_{\mathcal{M}}^{\text{community}}=\sum\limits_{m=1}^M C_m^{\text{grid}}\hspace{-1mm}+\hspace{-1mm}\sum\limits_{m=1}^M C_m^{\text{storage}}\hspace{-1mm}+\hspace{-1mm}\sum\limits_{m=1}^M\hspace{-1mm}C_m^{\text{operation}}\hspace{-1mm}-\hspace{-1mm}\sum\limits_{p=1}^P C_p^{\text{purchase}}\hspace{-1mm}. 
\end{equation}

We may formulate the whole community optimization problem that also includes sole energy consuming households as:

\vspace{2mm}
\hspace{-3.5mm}$\min\limits_{\{\mathbf{b}_m, \mathbf{r}_m, \mathbf{s}_m, \mathbf{a}_n \}_{m=1}^M,_{n=1}^N}\hspace{2mm} c_{\mathcal{M}}^{\text{community}}+$ \\

\hspace{19mm}$\sigma\hspace{-1mm}\sum\limits_{m=1}^M\hspace{-0.5mm}\sum\limits_{t=1}^T[b_m(t)\hspace{-0.5mm}+\hspace{-0.5mm}w_m(t)\hspace{-0.5mm}-\hspace{-0.5mm}u_m(t)\hspace{-0.5mm}-\hspace{-0.5mm}a_m(t)\hspace{-0.5mm}-\hspace{-0.5mm}r_m(t)]$,\\

\hspace{-3.5mm}such that the constraints in (\ref{eq:eq2})-(\ref{eq:eq8}) and (\ref{eq:eq14}) are satisfied, with the difference that (\ref{eq:eq14}) is reformulated for the whole community $\mathcal{N}$:
\vspace{-3mm}
\begin{equation}\label{eq:eq26}
\sum_{n=1}^N a_n(t)=0. 
\end{equation}
Just like in previous cases, we add a penalty term:\\ $\sigma\hspace{-1mm}\sum\limits_{m=1}^M\sum\limits_{t=1}^T[b_m(t)\hspace{-1mm}+\hspace{-1mm}w_m(t)\hspace{-1mm}-\hspace{-1mm}u_m(t)\hspace{-1mm}-\hspace{-1mm}a_m(t)\hspace{-1mm}-\hspace{-1mm}r_m(t)]$.

We denote by $c_p^{\text{indiv}}$ the cost paid by a household $p\in \mathcal{P}$ that participates in the community optimization. This cost would be equal to: $c_p^{\text{indiv}}=C_p^{\text{grid}}+C_p^{\text{operation}}+C_p^{\text{purchase}}$. For making sure that all sole energy consumers $p\in\mathcal{P}$ also benefit and reduce their energy cost the following constraint is imposed:
\vspace{-2mm}
\begin{equation}\label{eq:eq17}
c_p^{\text{indiv}}<\sum_{t=1}^Tu_p(t)\xi(t),\hspace{1mm} \forall \hspace{1mm} p\in \mathcal{P}.
\end{equation}
This constraint states that the cost paid by a sole energy consuming household that participates in the optimization should be less  than the cost that the household would pay by simply purchasing electricity directly from the utility company.

The final set of constraints for this problem is (\ref{eq:eq2})-(\ref{eq:eq8}), (\ref{eq:eq16}) and (\ref{eq:eq17}). The set of solution variables is $\{\mathbf{b}_m, \mathbf{r}_m, \mathbf{s}_m, \mathbf{a}_n \}_{m=1}^M,_{n=1}^N$.

It can be observed that the objective functions and the constraints of the formulated optimization problems possess linear relationships among the variables of the problems. Hence, we are able to model both optimization problems as linear programs. The solutions of the proposed optimization problems can be easily obtained through algorithms such as the interior point algorithm \cite{U21_2}.

\section{The coalitional game model for households owning RESs and ESSs}
\label{section:4}


The proposed coalitional cost optimization problems are using the energy storage systems and the renewable energy produced by the households $m\in\mathcal{M}$ in order to minimize the electricity cost. In this work we propose a characteristic coalitional game with transferable utility among the RESs and ESSs owning households.  

A coalitional game in characteristic form \cite {N14_2} is uniquely defined by the pair $(\mathcal{M},\upsilon)$, where $\mathcal{M}$ represents the set of players participating in the game, $|\mathcal{M}|=M$, and $\upsilon:2^{M}\to\mathbb{R}$ is the characteristic function of the game which quantifies the worth of a coalition. In the proposed optimization problem, the players $\mathcal{M}$ of the coalitional game are the households owning renewable energy production and storage systems. In order to minimize their electricity costs, these households may form different coalitional groups. We denote by $\mathcal{G}$, $\mathcal{G}\subseteq\mathcal{M}$, any non-empty subset of households from the set $\mathcal{M}$ that may form a coalitional group. The cardinality of the set $\mathcal{G}$ is $|\mathcal{G}|=G$. In case the coalition is formed by all the residences from the set $\mathcal{M}$, $\mathcal{G}=\mathcal{M}$, then this coalition is called the grand coalition. The worth of a coalition, which we generally denoted by $\upsilon(\mathcal{G})$, is a real value representing the total revenue received by the coalitional group for cooperating. 

In this work the worth (revenue) of a coalition is defined as the cost saved by a coalitional group in the collaborative  scenario in comparison with the total electricity cost that the members of that coalitional group would pay in the case of performing individual cost optimization. In this work we propose two coalitional optimization methods. Let us denote by $\upsilon(\mathcal{G})^{\text{ResEss}}$ the worth of a coalition for the case in which the optimization includes only the households $m\in\mathcal{M}$ owning RESs and ESSs. This worth is expressed as:  
\vspace{-2mm}
\begin{equation}\label{eq:eq18} 
\upsilon(\mathcal{G})^{\text{ResEss}}=\sum_{g=1}^Gc_g^{\text{indiv}}-c_{\mathcal{G}}^{\text{ResEss}}.
\end{equation}

Then, we denote by $\upsilon(\mathcal{G})^{\text{community}}$ the worth of a coalition for the case in which the optimization includes the whole community: RESs and ESSs owners and also the households $p\in\mathcal{P}$ that only consume electricity. In this case the worth is expressed as:
\vspace{-3mm}
\begin{equation}\label{eq:eq19} 
\upsilon(\mathcal{G})^{\text{community}}=\sum_{g=1}^Gc_g^{\text{indiv}}-c_{\mathcal{G}}^{\text{community}}.
\end{equation}
Here $c_g^{\text{indiv}}$ is computed using the formulation in (\ref{eq:eq21}), while $c_{\mathcal{G}}^{\text{ResEss}}$ and $c_{\mathcal{G}}^{\text{community}}$ are computed using (\ref{eq:eq22}) and (\ref{eq:eq23}), respectively.  The index of a household from the subset $\mathcal{G}$ that forms the coalitional group is denoted by $g$.

In coalitional games with transferable utility the worth of a coalition $\upsilon(\mathcal{G})$ is generally distributed among the members of the coalitional group using a fairness rule. There are various fairness methods in the literature for division of rewards like nucleolus, egalitarian, Shapley value \cite {N14_2,O15_2} for example. In this work we choose that the monetary revenues obtained by the coalition is divided among its members using the Shapley value \cite{O15_2}. The Shapley value is a method in which the worth of the coalition is distributed according to the amount of contribution that each player is bringing to the coalitional game and to the revenue. Shapley value does not depend on the order in which the players are joining the coalition. 

For a coalitional game defined by $(\mathcal{M},\upsilon)$, the Shapley value, $\Phi(\upsilon)$, assigns to each player \emph{$m\in\mathcal{M}$} a payoff $\Phi_m(\upsilon)$ given by the following expression:
\begin{equation}\label{eq:eq20} 
\Phi_m(\upsilon)\hspace{-0.5mm}=\hspace{-2mm}\sum_{\mathcal{G}\subseteq\mathcal{M}\setminus\{m\}}\hspace{-2mm}\frac{G!(M-G-1)!}{M!}[\upsilon(\mathcal{G}\cup\{m\})-\upsilon(\mathcal{G})],
\end{equation}
where the sum is computed over all possible subsets $\mathcal{G}$ (even of single players) of $\mathcal{M}$ not containing player $m$. The payoff of a household taking part in a coalition represents the fraction of the total revenue of that coalition, that is achieved through the participation of that household in the coalitional game. In the proposed game all households from set $\mathcal{M}$ participate in the coalition. The distributed payyofs must obey the following rule:$\sum\limits_{m\in\mathcal{M}}\Phi_m(\upsilon)=\upsilon(\mathcal{M})$.

\section{Simulation Examples}
\label{section:5}

\begin{figure*}[!t]

	\centering
		\includegraphics[width=.8\textwidth]{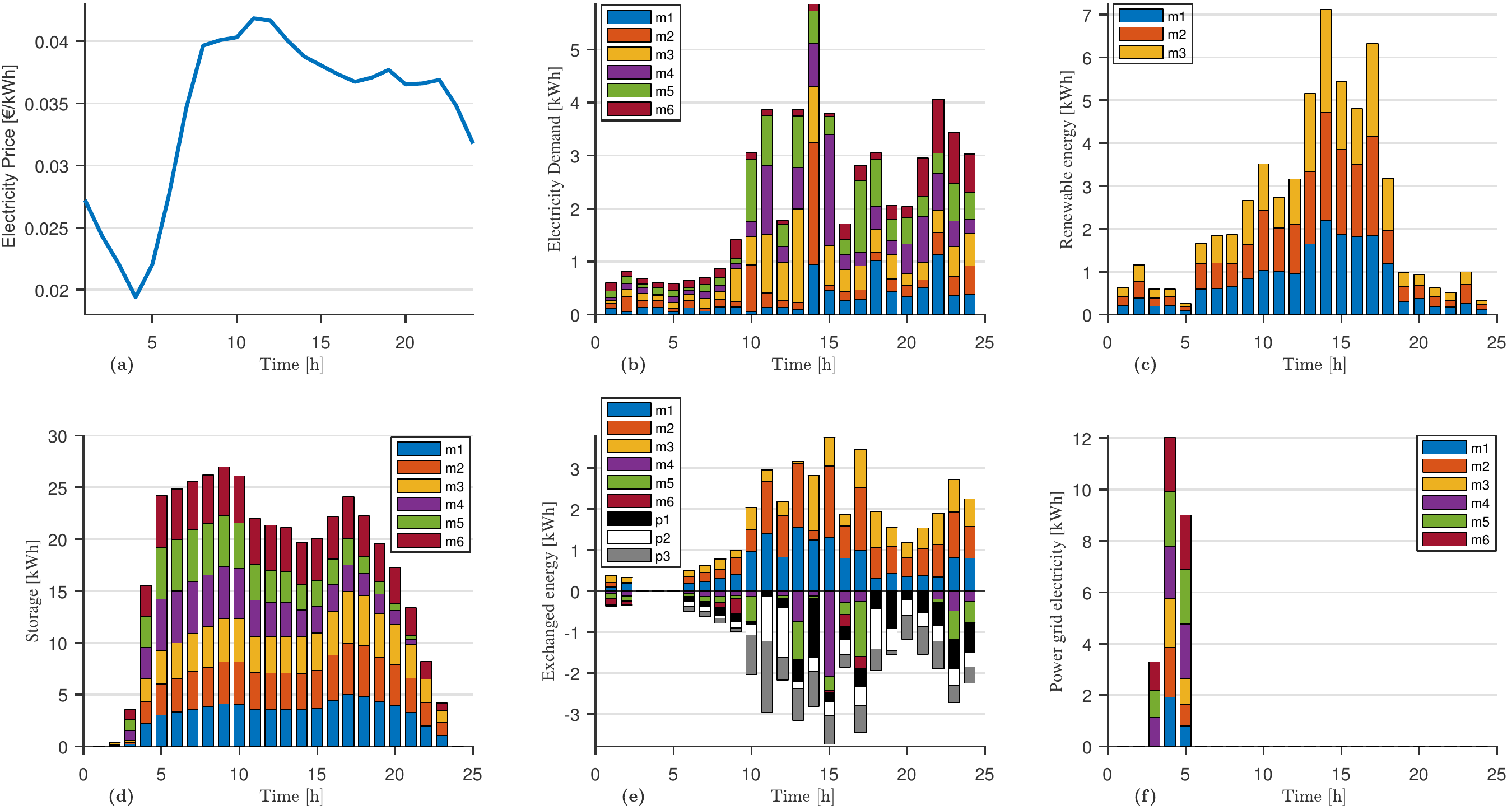}
		\caption{(a) 24-hours utility company electricity price. (b) 24-hours electricity demand of the households owning RESs and ESSs. (c) 24-hours renewable energy production of households owning RESs. (d) 24-hours ESSs profiles. (e) 24-hours amounts of electricity exchange between all households in the community. The positive blocks show amounts of energy provided by some households to the others, whereas the negative blocks show amounts of energy received by the other households. (f) Electricity purchased from the utility company by RESs and ESSs owners. Households buy energy from the utility during the hours when the price is low.}\label{fig:Solutions1}
\end{figure*} 
In this section we present simulation examples and quantitative results that demonstrate the performance and cost savings achieved by the proposed methods. For simulating and testing the performance of the proposed methods we considered a smart grid community composed of \emph{$N$}= 9 households. The set $\mathcal{M}$ of residences is composed of 6 households. In our simulations we denote these households by $\mathcal{M}=\{m1, m1,\ldots\,m6\}$. Out of this set, 3 households, $\{m1, m2, m3\}$, own RESs as well as ESSs, while the other 3 households, $\{m4, m5, m6\}$, own ESSs only.  The remaining \emph{$|\mathcal{P}|$}=3 households in the community are sole energy consumers: $\mathcal{P}=\{p1, p2, p3\}$. For the presented results, except those in Fig. 4, we considered equal $ESS_{\emph{m}}$ capacities, $C_m$, of 5kWh. We also assumed equal charging/discharging rates, $\rho_m$, of 2kWh. The storage loss factor is assumed to be $\eta_m=0.001, m=1,\ldots,M$. We perform simulations over a time horizon $\mathcal{T}$ of length \emph{T}=24 hours divided into hourly time slots. The utility company pricing data, $\boldsymbol{\xi}$, is true pricing data taken from Finnish Nord Pool Spot database \cite{P16_2} for May 2013. Other pricing data used in simulations is the following: $\pi$=0.0001 $\euro$ per each kWh of charged/discharged energy, $\sigma$=0.001 $\euro$ per each kWh of renewable energy not stored or consumed, $\tau$=0.0001 $\euro$ per each kWh of electricity transferred or received by each household. The price of energy sold by the RESs and ESSs owners to the sole energy consumers is $\boldsymbol{\lambda}=\alpha\boldsymbol{\xi}, \alpha=0.9$, i.e. the sole consumers get a 10\% discount compared to the price offered by the utility company. The influence of $\alpha$ over the performances of the proposed methods is depicted in Fig. 3. A high variation of the chosen values for the prices mentioned above ,$\pi,\rho,\tau$, may have an influence of about $\pm1.2\%$ on the cost savings resulted from the proposed methods. 
 
In this paper we show simulation results of the proposed methods over 31 days. In these simulations we assumed empty storages, $s_m(0)$=0, $m=1,\ldots,M$, at the beginning of the first 24-hours optimization period. Further, each proposed optimization method updated its corresponding initial storage values, for each new optimization period, with the storage levels resulted at the end of the previous period $s_m(0)=s_m(24)$. Hence, the proposed optimization methods may have had different initial storage levels at the beginning of the 24-hours optimization periods from the 31-days simulation. In order to simulate the 24-hours electricity demands of the households, $\mathbf{u}_n$, we used the load modeling framework presented in \cite{Q17_2}. We assumed residences with the following numbers of inhabitants: \{3, 4, 4, 2, 5, 4, 3, 4, 2\}. We  considered that the households were equipped with wind energy producing systems. For simulating renewable energy data values we used the following mathematical model to approximate the power generated by a wind turbine \cite{R18_2}: $P_w=(1/2)DK_pAV^3$, where $D$ represents the air density, $K_p$=0.3 is the power coefficient of the turbine, $A$ is the swept area of the turbine, $A=3.14R^2,R=2.63$m, and $V$ is the wind speed. For this, we used true weather data for May 2013 in Helsinki region \cite{S19_2}. To represent the energy production of each household, we added a small variation to the values calculated with the above formulation. For solving the linear programs resulted from modeling the proposed optimization problems we used the CVX package for convex optimization \cite{T20_2}. The results of the performed simulations are presented further.

Fig. 1 (a)-(f) shows an example of 24-hours input data and results of the \textit{whole community optimization problem}. We used the pricing and weather data for May 17 \cite{P16_2,S19_2}. Fig. 1.(a) shows the electricity prices, $\boldsymbol{\xi}$, in $\euro/kWh$. Fig. 1.(b) shows the individual and cumulated hourly electricity demands of those households that own RESs and ESSs: $\{m1,\ldots, m6\}$. Fig. 1.(c) shows the individual and cumulated hourly renewable energy production of households $\{m1,m2, m3\}$ that own RESs. Fig. 1.(d) shows the individual and cumulated hourly profiles of the ESSs, i.e how much energy exists in the energy storages of households $m\in \mathcal{M}$ at the end of each hour. Fig. 1.(e) shows the amounts of energy exchanged during each hour by all the households in the community. At each time-slot, the positive blocks show the amounts of energy that are provided by some households to the others, whereas the negative blocks show the amounts of energy that are being received by the other households. In this example it can be observed that in most cases the households which produce renewable energy provide energy to those households that own ESSs only and to the sole energy consuming households in the network. Fig. 1.(f) shows the amounts of energy purchased from the power grid by those households, $m\in \mathcal{M}$, that own RESs and ESSs. It can be seen that the households buy energy from the utility company during those hours when the electricity price is low.      
\begin{figure}[!t]
	\centering
		\includegraphics[width=\linewidth]{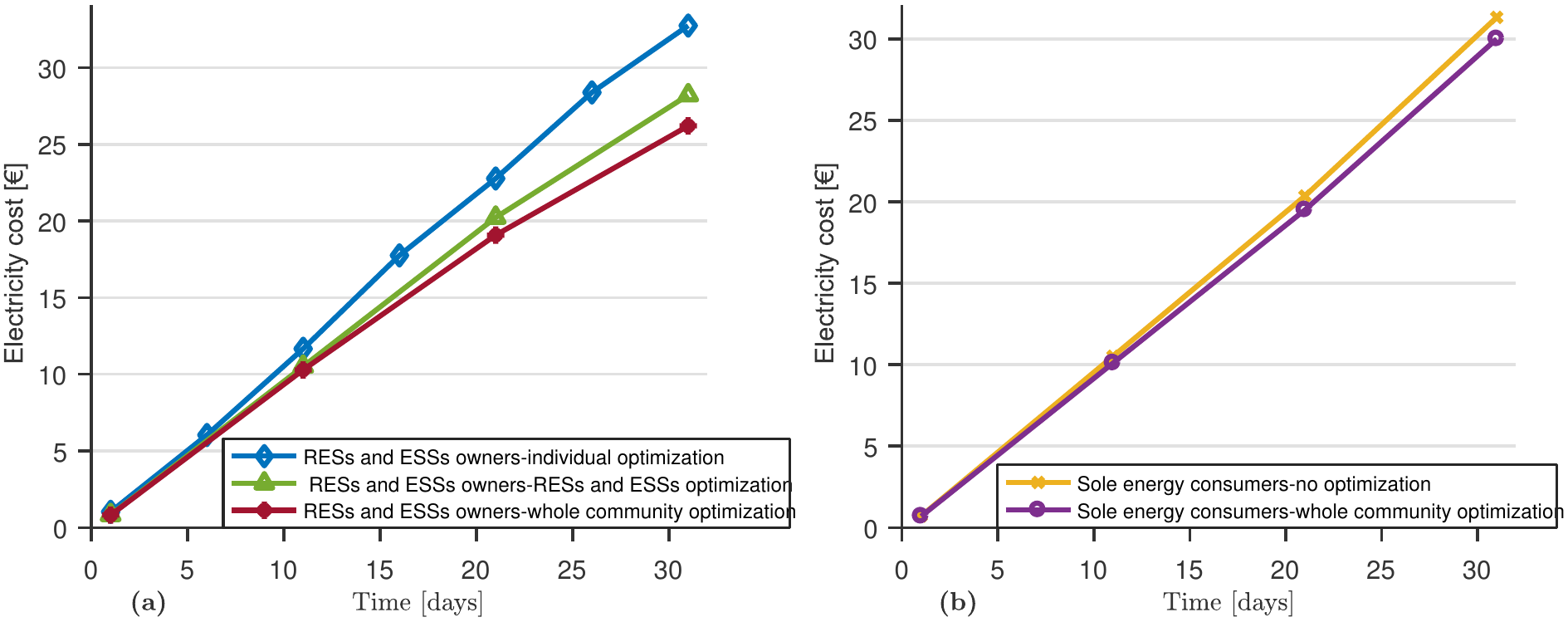}
		\caption{Comparison between the 31-days cumulative electricity costs of the RESs and ESSs owners (a) and of the sole energy consumers (b), obtained through the proposed optimization methods. The coalitional cost minimization for the households owning RESs and ESSs results in a cost reduction of about 13.5\% in comparison to the individual cost optimization. In case of whole community coalitional optimization, the RESs and ESSs owners obtain a cost reduction of about 20\% compared to the individual optimization. The sole energy consummers obtain a cost reduction of about 4\%.}\label{fig:Solutions2}
\end{figure} 

Fig. 2 shows in (a) the 31-days cumulative electricity costs of the households $m \in \mathcal{M}$ that own RESs and ESSs  and in (b) the cumulative electricity costs of the households $p\in \mathcal{P}$ that are sole energy consumers.  It can be observed in (a) that the \textit{ coalitional cost minimization method for the households owning RESs and ESSs} (28.21$\euro{}$) provides for these residences a cost reduction of about 13.5\% versus the \textit{individual cost optimization} (32.74$\euro{}$). The \textit{whole community cost minimization} problem (26.2$\euro{}$) provides an additional 6.5\% reduction in cost, resulting in an overall 20\% cost reduction for the RESs and ESSs owners. In (b) it can be observed that the sole energy consumers, through participating to the \textit{whole community cost optimization} (30.03$\euro{}$) reduce their costs, in this example by approximately 4\% in comparison to buying all electricity from the utility company only (31.32$\euro{}$). 
\begin{figure}[!t]
	\centering
		\includegraphics[width=.7\linewidth]{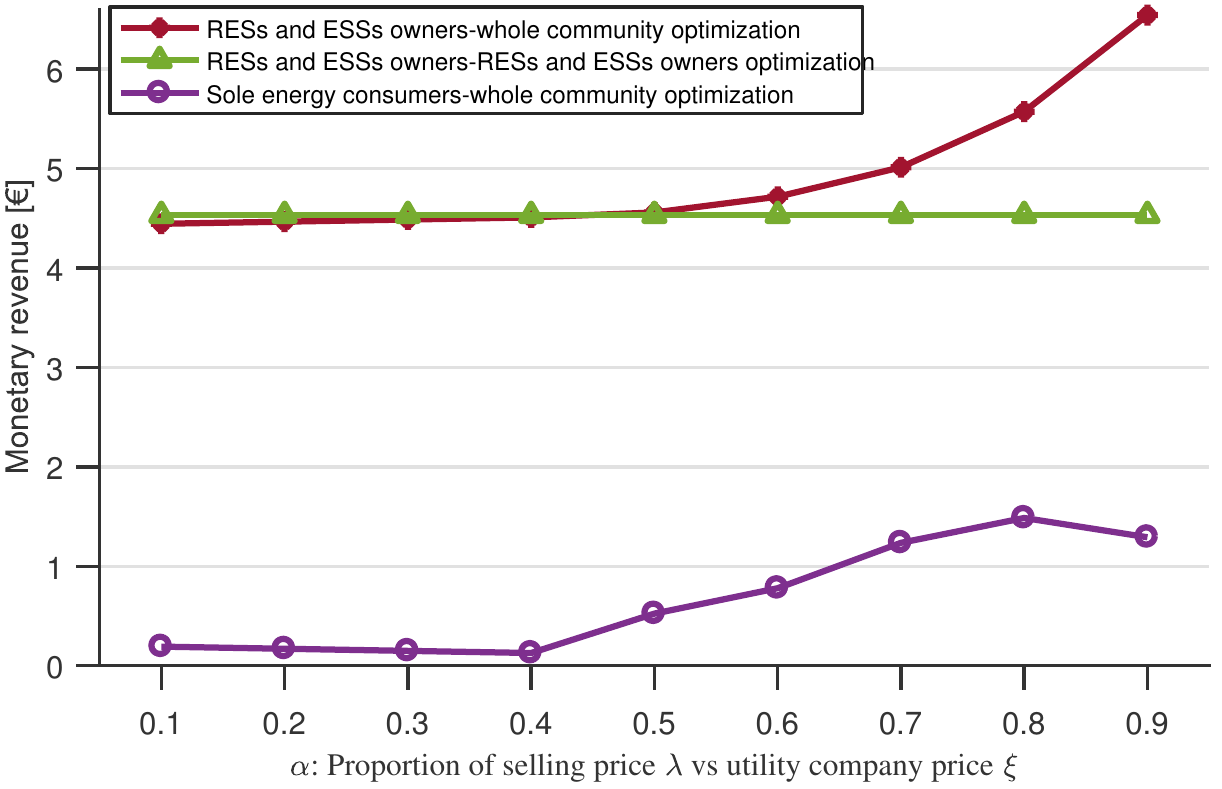}
		\caption{The variation of the 31-days cumulated monetary revenues of the proposed coalitional problems in comparioson to $\alpha$, the proportion of the price of electricity sold by the RESs and ESSs owners $\boldsymbol{\lambda}$ versus the utility company price $\boldsymbol{\xi}$. The whole community optimization method  obtains positive results for $\alpha\geq 0.4$. The highest gain is obtained by the sole energy consumers when $\alpha=0.8$ and by the RESs and ESSs owners when $\alpha=0.9$.}\label{fig:Solutions3}
\end{figure} 

In Fig. 3 we want to show how the 31-days cumulated monetary gains of the proposed coalitional optimization problems are varying for different values of the set of prices $\boldsymbol{\lambda}$. As specified in Section II, $\boldsymbol{\lambda}$ is the set of electricity prices corresponding to the energy sold by the RESs and ESSs owners to the sole energy consuming households. The values of $\boldsymbol{\lambda}$ are smaller than the utility company set of electricity prices, $\boldsymbol{\xi}$: $\boldsymbol{\lambda}=\alpha\boldsymbol{\xi}$, where $\alpha$: $0\leq\alpha\leq 1$. In this simulation we used the following values for $\alpha$: \{0.1,               \ldots, 0.9\}. It can be observed from Fig. 3 that the \textit{whole community optimization} obtains positive results for $\alpha\geq 0.4$. The highest cost reduction, of about 20\% from the \textit{individual optimization}, is obtained  by the RESs and ESSs owners when the price of sold energy is $\boldsymbol{\lambda}= 0.9\boldsymbol{\xi}$, i.e $\alpha$=0.9. For the sole energy consumers the highest cost reduction, of about 5\% from from cost of purchase from utility company, is obtained when $\alpha=0.8$. In the \textit{whole community cost optimization} the cost reduction of the sole energy consumers is a consequence of the cost minimization performed by the RESs and ESSs owners. Households $m \in \mathcal{M}$ minimize their cost by selling energy to households $p \in \mathcal{P}$. The higher are the prices $\boldsymbol{\lambda}$, the higher are the amounts of energy sold by the households $m\in\mathcal{M}$. Hence, the cost reduction of the sole energy consumers gets higher with the amounts of energy purchased from RESs and ESSs owners. 
\begin{figure}[!t]
	\centering
		\includegraphics[width=.7\linewidth]{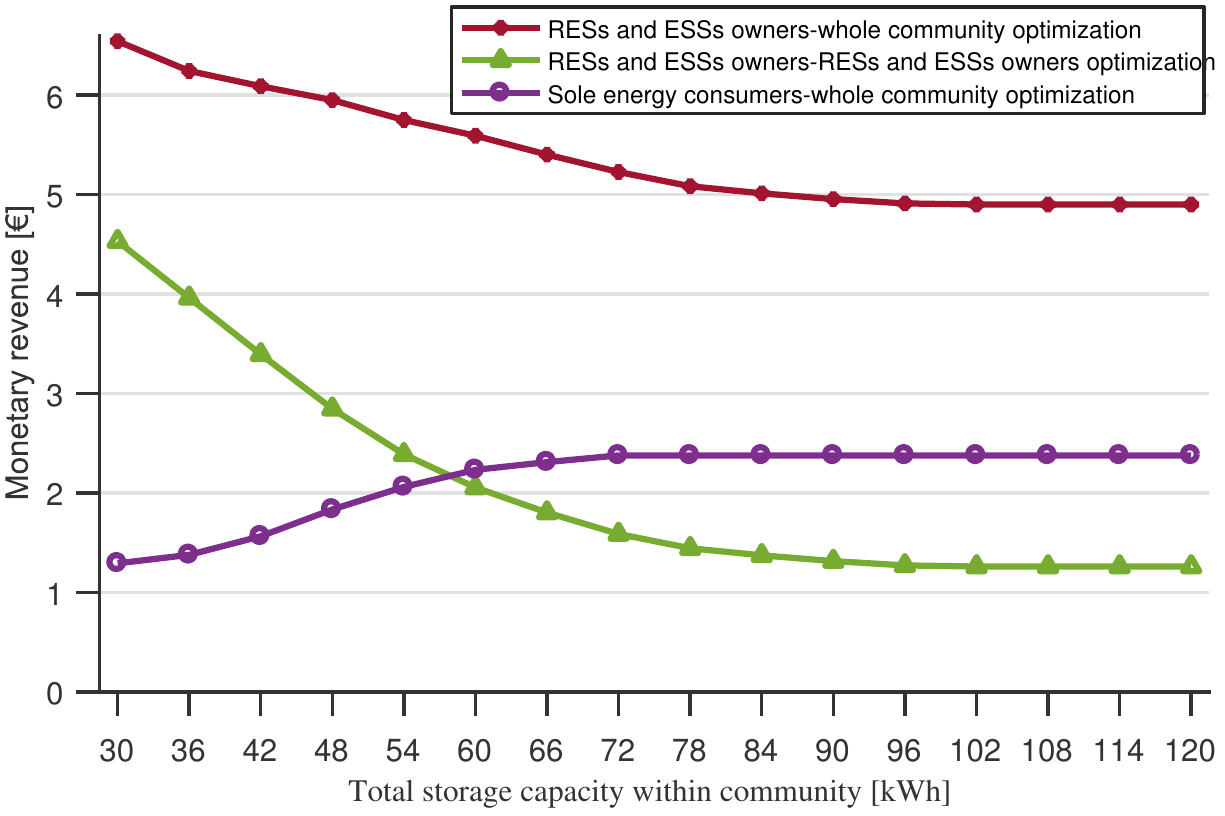}
		\caption{The variation of the 31-days cumulated monetary revenues obtained by the proposed coalitional methods versus the storage capacity existing whithin the community. The monetary gains of the RESs and ESSs owners decrease with the increase of the storage capacity for both proposed coalitional cost minimization methods. The gains of the sole energy consumers increases with the increase of storage capacity.}\label{fig:Solutions5}
\end{figure} 

Fig. 4 shows the variation of the 31-days cumulated the monetary gains of the proposed coalitional optimization methods
versus the storage capacity existing within the community. Initially we assumed capacities of 5kWh for each $ESS_m$, resulting in an overall storage capacity of 30kWh within the community. Then we increased the storage capacity of each $ESS_m$ by 1kWh, until reaching 20kWh, which resulted in an overall 120kWh storage capacity within the community. It can be observed that the monetary gains of the RESs and ESSs owners are actually decreasing while the capacity increases. This is due of the fact that the result of the \textit{individual cost optimization} is used as a comparison benchmark in (\ref{eq:eq18}) and (\ref{eq:eq19}) for the calculation of the revenues of the coalitional optimizations. The increase of the ESSs capacities allows for better performance of the \textit{individual cost minimization} problem. Hence, the coalitional optimizations' gains for RESs and ESSs owners are not so significant anymore. However, we can see that the gains of the proposed coalitional methods don't totally decrease. Actually, even with very large storage capacities, the proposed coalitional methods still obtain significant monetary gains in comparison to the individual optimization. Contrary to the gains of RESs and ESSs owners, it can be observed that the gains of the sole energy consumers become larger as the storing capacity within the community increases.
 
The bar plot in Fig. 5 depicts the daily cost savings of the community and of the individual households for the case in which the cooperation is done according to the proposed \textit{whole community coalitional game optimization}. Fig. 5(a) shows the cost savings of the RESs and ESSs owners. The daily cost savings are divided among the members of the community as shown by the blocks composing each bar. Each block indicates the amount of monetary payoff received by a household as resulted the Shapley value calculation. We can observe that a big part of the daily cost savings are allocated to the renewable energy producers. It can also be observed that during five days of the month, days 3,6,8,18 and 23 the coalitional optimization does not perform better than the individual optimization. This is also due to the fact that, as mentioned before, the storages of the EESs owners do not have same initial values for the individual optimization as for the coalitional optimization. During some days, the initial storage values may be full for the individual optimization, while for the coalitional optimization the initial storage values may be empty. For the majority of the days of the month each household in the community receives some monetary revenue for their participation in the coalitional optimization. Fig. 5(b) shows the cost savings of the sole energy consuming households as resulted from the optimization, by buying energy from the RESs and ESSs owners at a lower price than the one offered by the utility company. Each of the sole energy consuming household obtains daily cost savings.  
\begin{figure}[!t]
	\centering
		\includegraphics[width=.8\linewidth]{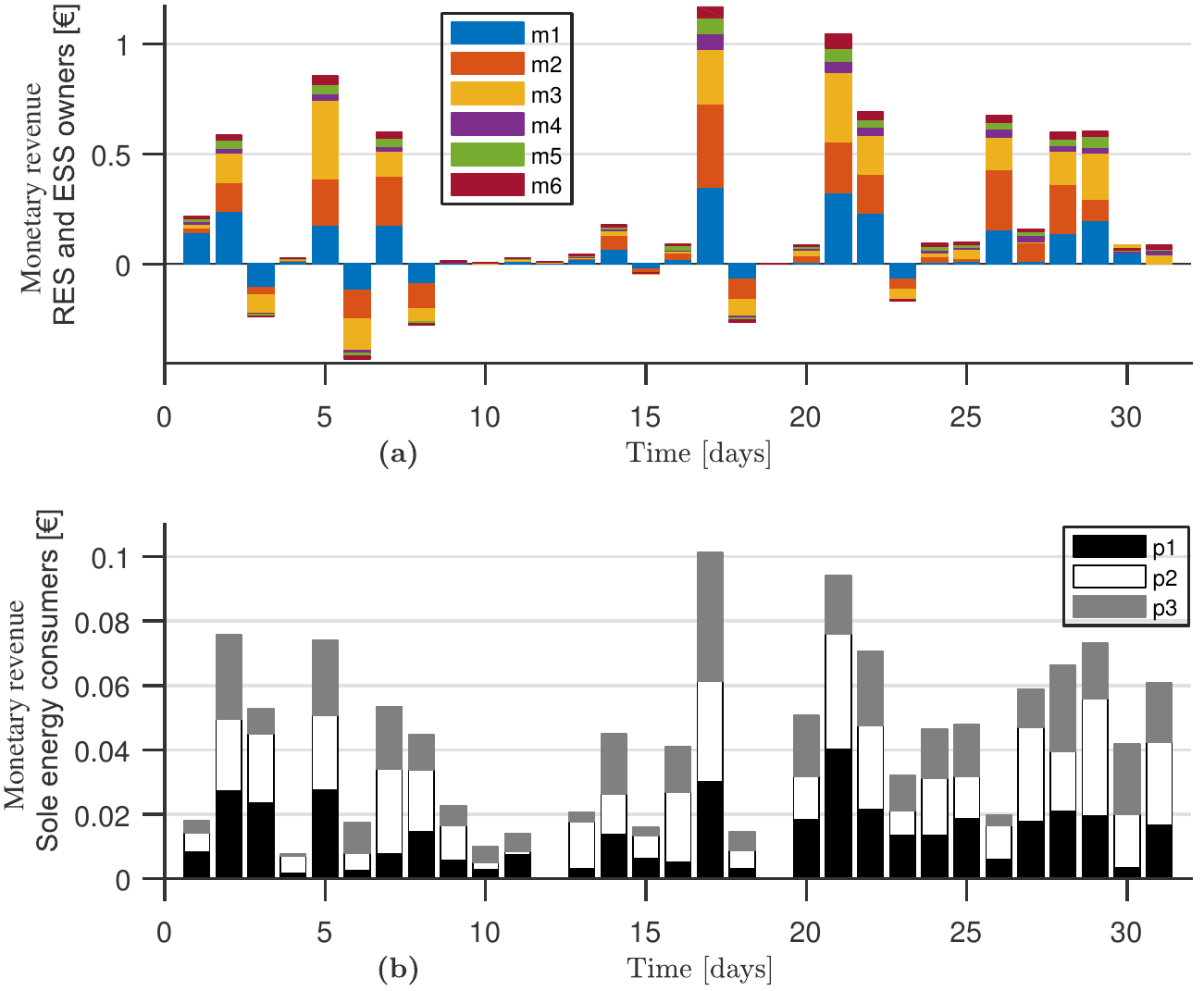}
		\caption{(a) The daily monetary revenues achieved by the RESs and ESSs owning households by cooperating in the whole community optimization. The blocks composing a bar show the payoffs received the households from the coalition as distributed by the Shapley value method. (b) The daily monetary revenues achieved by the sole energy consumming households through participating in the whole community optimization.}\label{fig:Solutions5}
\end{figure}

\section{Conclusions}
\label{section:6}
In this paper we proposed two novel coalitional game theory based optimization methods for minimizing the cost of electricity consumed by households in a smart grid community. We first proposed a coalitional optimization through which the RESs and ESSs owning households from the community can form coalitions for exchanging energy and sharing their renewable energy and storage units. Simulation examples showed that through participating to this  coalitional optimization game, these households can reduce their energy consumption costs by roughly 13.5\%. We then proposed a secondary coalitional cost optimization model through which the RESs and ESSs share their renewable energy and storage units among themselves, but also sell electricity to the pure consuming households in the community. We show that by selling electricity at a cost lower than the one offered by the utility company, the RESs and ESSs owners may further reduce their costs with about 6.5\%, resulting in a total 20\% cost reduction in comparison with the cost incurred for individually optimizing their cost. Also, by participating to the optimization, the sole energy consuming households may also obtain electricity cost reductions up to 5\%.

\bibliographystyle{IEEEtran}

\bibliography{mybib_2NDJournal} 

\begin{thebibliography}{10}
\providecommand{\url}[1]{#1}
\csname url@samestyle\endcsname
\providecommand{\newblock}{\relax}
\providecommand{\bibinfo}[2]{#2}
\providecommand{\BIBentrySTDinterwordspacing}{\spaceskip=0pt\relax}
\providecommand{\BIBentryALTinterwordstretchfactor}{4}
\providecommand{\BIBentryALTinterwordspacing}{\spaceskip=\fontdimen2\font plus
\BIBentryALTinterwordstretchfactor\fontdimen3\font minus
  \fontdimen4\font\relax}
\providecommand{\BIBforeignlanguage}[2]{{%
\expandafter\ifx\csname l@#1\endcsname\relax
\typeout{** WARNING: IEEEtran.bst: No hyphenation pattern has been}%
\typeout{** loaded for the language `#1'. Using the pattern for}%
\typeout{** the default language instead.}%
\else
\language=\csname l@#1\endcsname
\fi
#2}}
\providecommand{\BIBdecl}{\relax}
\BIBdecl

\bibitem{A1_2}
H.~Wang and J.~Huang, ``Cooperative planning of renewable generations for
  interconnected microgrids,'' \emph{IEEE Trans. on Smart Grid}, vol.~7, no.~5,
  pp. 2486--2496, Sept 2016.

\bibitem{B2_2}
------, ``Incentivizing energy trading for interconnected microgrids,''
  \emph{available as early access article in IEEE Trans. on Smart Grid}, 2016.

\bibitem{C3_2}
K.~Rahbar, C.~C. Chai, and R.~Zhang, ``Energy cooperation optimization in
  microgrids with renewable energy integration,'' \emph{available as early
  access article in IEEE Trans. on Smart Grid}, 2016.

\bibitem{D4_2}
A.~Ouammi, H.~Dagdougui, L.~Dessaint, and R.~Sacile, ``Coordinated model
  predictive-based power flows control in a cooperative network of smart
  microgrids,'' \emph{IEEE Trans. on Smart Grid}, vol.~6, no.~5, pp.
  2233--2244, Sept 2015.

\bibitem{E5_2}
Y.~Wang, S.~Mao, and R.~M. Nelms, ``On hierarchical power scheduling for the
  macrogrid and cooperative microgrids,'' \emph{IEEE Trans. on Industrial
  Informatics}, vol.~11, no.~6, pp. 1574--1584, Dec 2015.

\bibitem{F6_2}
L.~Xiao, N.~B. Mandayam, and H.~V. Poor, ``Prospect theoretic analysis of
  energy exchange among microgrids,'' \emph{IEEE Trans. on Smart Grid}, vol.~6,
  no.~1, pp. 63--72, Jan 2015.

\bibitem{G7_2}
D.~Gregoratti and J.~Matamoros, ``Distributed energy trading: The
  multiple-microgrid case,'' \emph{IEEE Trans. on Industrial Electronics},
  vol.~62, no.~4, pp. 2551--2559, April 2015.

\bibitem{H8_2}
Y.~Wang, W.~Saad, Z.~Han, H.~V. Poor, and T.~Başar, ``A game-theoretic
  approach to energy trading in the smart grid,'' \emph{IEEE Trans. on Smart
  Grid}, vol.~5, no.~3, pp. 1439--1450, May 2014.

\bibitem{L12_2}
J.~Rajasekharan and V.~Koivunen, ``Cooperative game-theoretic approach to load
  balancing in smart grids with community energy storage,'' in \emph{2015 23rd
  European Signal Processing Conference (EUSIPCO)}, Aug 2015, pp. 1955--1959.

\bibitem{I9_2}
I.~Atzeni, L.~G. Ordóñez, G.~Scutari, D.~P. Palomar, and J.~R. Fonollosa,
  ``Noncooperative and cooperative optimization of distributed energy
  generation and storage in the demand-side of the smart grid,'' \emph{IEEE
  Trans. on Signal Processing}, vol.~61, no.~10, pp. 2454--2472, May 2013.

\bibitem{K11_2}
W.~Saad, Z.~Han, and H.~V. Poor, ``Coalitional game theory for cooperative
  micro-grid distribution networks,'' in \emph{2011 IEEE International
  Conference on Communications Workshops (ICC)}, 2011, pp. 1--5.

\bibitem{J10_2}
W.~Tushar, C.~Yuen, D.~B. Smith, N.~U. Hassan, and H.~V. Poor, ``A canonical
  coalitional game theoretic approach for energy management for nanogrids,'' in
  \emph{2015 IEEE Innovative Smart Grid Technologies - Asia (ISGT ASIA)}, Nov
  2015.

\bibitem{M13_2}
A.~Chiş, J.~Lundén, and V.~Koivunen, ``Optimization of plug-in electric
  vehicle charging with forecasted price,'' in \emph{2017 IEEE International
  Conference on Acoustics, Speech and Signal Processing (ICASSP)}, March 2017.

\bibitem{U21_2}
Y.~Ye, \emph{Interior point algorithms: theory and analysis}.\hskip 1em plus
  0.5em minus 0.4em\relax Wiley, 1997.

\bibitem{N14_2}
M.~Osborne and A.~Rubinstein, \emph{A Course in Game Theory}.\hskip 1em plus
  0.5em minus 0.4em\relax MIT Press, 1994.

\bibitem{O15_2}
W.~Saad, Z.~Han, M.~Debbah, A.~Hjorungnes, and T.~Basar, ``Coalitional game
  theory for communication networks,'' \emph{IEEE Signal Processing Magazine},
  vol.~26, no.~5, pp. 77--97, September 2009.

\bibitem{P16_2}
O.~A. www.nordpoolspot.com.

\bibitem{Q17_2}
I.~Richardson, M.~Thomson, D.~Infield, and C.~Clifford, ``Domestic electricity
  use: A high-resolution energy demand model,'' \emph{Energy and Buildings},
  vol.~42, no.~10, p. 1878–1887, 2010.

\bibitem{R18_2}
T.~Burton, D.~Sharpe, N.~Jenjkins, and E.~Bossanyi, \emph{Wind Energy
  Handbook}.\hskip 1em plus 0.5em minus 0.4em\relax Wiley, 2001.

\bibitem{S19_2}
O.~A. www.ilmatieteenlaitos.fi.

\bibitem{T20_2}
M.~Grant and S.~Boyd.\hskip 1em plus 0.5em minus 0.4em\relax CVX: Matlab
  Software for Disciplined Convex Programming, Version 2.0 Beta, Sept 2013.

\end{thebibliography}

\end{document}